\documentclass[12pt]{article}
\usepackage{fullpage}
\usepackage{graphicx}
\usepackage[breaklinks=true]{hyperref}
\usepackage{xcolor}
\hypersetup{
    colorlinks,
    linkcolor={red!50!black},
    citecolor={blue!50!black},
    urlcolor={blue!80!black}
}

\usepackage{amsmath}
\usepackage{amssymb}
\usepackage{mathtools}
\usepackage{amsthm}
\usepackage[capitalize,noabbrev]{cleveref}

\usepackage{natbib}
\usepackage[english]{babel}
\usepackage{multirow}

\newcommand{\vecc}{\mathbf{c}}
\newcommand{\x}{\mathbf{x}}

\newcommand{\vecY}{\mathbf{y}}
\newcommand{\vecYu}{\bar{\mathbf{y}}}
\newcommand{\Id}{I}

\newcommand{\Var}[1]{\mathrm{Var}\left[{#1}\right]}
\newcommand{\Esp}[1]{\mathrm{E}\left[{#1}\right]}

\newcommand{\Normal}[1]{\mathrm{N}\left(#1 \right)}

\newcommand{\mn}{m_N} 
\newcommand{\vn}{s_N^2} 
\newcommand{\CN}{\mathbf{C}_N} 
\newcommand{\Cn}{\mathbf{C}_n} 
\newcommand{\An}{\mathbf{A}_n} 
\newcommand{\LaN}{\boldsymbol{\Lambda}_N} 
\newcommand{\Lan}{\boldsymbol{\Lambda}_n} 
\newcommand{\pvar}{\sigma^2} 
\newcommand{\xu}{\bar{\mathbf{x}}} 
\newcommand{\Xset}{\mathbb{X}} 

\title{Gearing Gaussian process modeling and sequential design towards noisy simulators}
\author{Micka\"{e}l Binois\thanks{Inria, Université Côte d’Azur, CNRS, LJAD, Sophia Antipolis, France \href{mailto:mickael.binois@inria.fr}{\tt mickael.binois@inria.fr}}
\and
Arindam Fadikar\thanks{Argonne National Laboratory, Decision and Infrastructure Sciences, U.S.A.}
\and Abby Stevens\footnotemark[2]
}
\date{}

\begin{document}

\maketitle

\abstract{
This chapter presents specific aspects of Gaussian process modeling in the presence of complex noise. Starting from the standard homoscedastic model, various generalizations from the literature are presented: input varying noise variance, non-Gaussian noise, or quantile modeling. These approaches are compared in terms of goal, data availability and inference procedure. A distinction is made between methods depending on their handling of repeated observations at the same location, also called replication. 
The chapter concludes with the corresponding adaptations of the sequential design procedures. These are illustrated in an example from epidemiology.
}

\section{Introduction}

Accurately reproducing real-world dynamics often requires stochastic simulators,
particularly in fields like epidemiology, operations research, and
hyperparameter tuning. In these contexts it becomes important to distinguish between aleatoric uncertainty -- arising from noise in observations, from epistemic
uncertainty -- stemming from uncertainty in the model. The former is sometimes called intrinsic uncertainty while the latter is referred to as extrinsic
uncertainty, see e.g., \cite{Ankenman2010}.

Gaussian process (GP) based surrogate methods (see, e.g., \cite
{Rasmussen2006,Gramacy2020}) can be easily adapted from deterministic to noisy settings while maintaining strong predictive power, computational efficiency, and analytical tractability. Even in the deterministic setup, it is common to add a small diagonal \emph{nugget} (also known as a \emph{jitter}) term to the covariance matrix of the GP
equations to ease its numerical inversion. It is also interpreted as a
regularization term, especially in the reproducing kernel Hilbert space (RKHS)
context, see, e.g., \cite{kanagawa2018gaussian}. This can be contrasted to
the use of pseudo-inverses, which reverts to interpolation, see for instance the discussion by \cite
{mohammadi2016analytic}. Here we will prefer the term \emph{noise variance} to
relate it to intrinsic uncertainty, and also because the nugget effect has a
different meaning in the kriging literature (see e.g., \cite{Roustant2012}).

Surrogates are often developed to replace calls to expensive simulators, and are hence used in a variety of contexts: uncertainty propagation, sensitivity analysis, Bayesian
optimization, calibration, or inversion, see e.g., \cite{Gramacy2020}. Typically, from an initial set of observations, additional ones are added sequentially based on a suitable criterion. While relying on deterministic
tools for these is possible, ignoring the noise variance may be
detrimental to understanding complex dynamics and further adaptations are generally needed. For example, in both modeling
and sequential design contexts, experiments are often replicated (i.e., different outputs are repeatedly generated from the same inputs due to the stochasticity). If
this should be prevented in the noiseless context, it is an efficient tool to
learn the noise features. Savings in total evaluation budget and
computational speed may be obtained by adaptively balancing replication in
the input space.

The goal of this chapter is to provide detailed treatments for situations when the
intrinsic uncertainty becomes significant, resulting in low signal-to-noise ratios. We address both the modeling aspects and then the sequential design strategies under these conditions. In
particular, we focus on alternatives for modeling the corresponding noise
process, with various degrees of complexity and depending on data
availability, with or without replication. While detailing only a few of the possible models and
associated references, we refer for instance to \cite
{kleijnen2015design,Gramacy2020,baker2022analyzing} for complementary and
more exhaustive viewpoints on the options. Sequential design procedures to enrich an initial data set are detailed, again taking into account the option to replicate.

This chapter is organized as follows: Section \ref{sec:noisygps} details the
simple adaptations for noisy modeling, disambiguating different contexts
while keeping analytical tractability. Section \ref{sec:advnoisy} discusses
further refinements requiring approximations. Section \ref
{sec:seqdesign} details some modifications towards sequential designs.
Finally, Section \ref{sec:concl} summarizes the main concepts and opens to
remaining challenges.

\subsection{Related frameworks}

For numerical simulators, the aleatoric behavior is generally controlled by a
random seed argument. Given this seed information, seen as a categorical variable,
the simulator becomes deterministic. As a result, dedicated covariance
functions can be defined, see for instance \cite{Chen2012a,Pearce2022} for
common random numbers (CRN) oriented ones, or \cite
{Zhou2011a,roustant2020group,zhang2020latent,Deshwal2021} for the more
general options.

A similar context is when additional
environmental variables are present: they correspond to variables that may
not be controlled but whose values are known. Given the control and
environmental variables, the output value remains deterministic. 
An additional related instance is when realizations of the random noise
process are available, which may be used to estimate the behavior of the
environmental variables, see e.g., \cite{el2023feasible}. This pathwise or
trajectory interpretation is further explored and exploited, e.g., by \cite
{Wilson2021,fadikar2023trajectory}. 

Considering input noise is an alternative to replace the environmental
variables effect, with a given probability distribution on the inputs. For
further analysis, a robust GP is defined which integrates out the
environmental variables or input noise. Without further assumptions, Monte
Carlo integration is required to do so, while analytical expressions are
available for simpler setups, typically with Gaussian noise, see e.g., \cite
{girard2002gaussian,McHutchon2011,janusevskis2013simultaneous,Qing2022}. More
generally, GPs may also be defined directly on distributions, see e.g., \cite
{bachoc2020gaussian}.

Notice that these differ from the i.i.d.\ noise assumption that we make with
noisy GPs, with no extra information to exploit (irreducible noise). These
setups, however, may be combined, e.g., having both input and output noise.

\section{Noisy GP modeling}
\label{sec:noisygps}

While the modification to account for noise in Gaussian processes is relatively straightforward, inference for an unknown input-varying noise variance function comes with some challenges. Options are discussed and illustrated on a data set.

\subsection{Extending the constant noise framework}

Denoting $Y$ the GP model of the simulator, recall the classical GP predictive equations for $N$ observations with i.i.d.\
additive Gaussian noise, $Y_i(\x_i) = f(\x_i) + \varepsilon_i$, $1 \leq i \leq N$, where
$\varepsilon_i \sim \Normal{0, r(\x_i)}$ with $r(\cdot)$ is the
(positive) variance function:
\begin{align}
\mn(\x) &= \Esp{Y(\x)|y_1, \dots, y_N} = \vecc(\x)^\top (\CN + \LaN)^{-1} \vecY \\
\vn(\x) &= \Var{Y(\x)|y_1, \dots, y_N} = r(\x) + \pvar \left(c(\x, \x) - \vecc(\x)^\top (\CN
 + \LaN)^{-1} \vecc(\x) \right)
 \label{eq:hetGP}
\end{align}
with $\LaN = \sigma^{-2} Diag\left(r(\x_1), \dots, r(\x_N) \right)$.

Compared to the standard homoscedastic GP model, with constant noise, that is
$\LaN = \nu \Id$, here the noise variance can vary with $\x$ while keeping the
analytical tractability of the GP prediction. Note that the predictive variance includes $r
(\x)$, removing which corresponds to $\Var{f(\x)|y_1, \dots, y_N}$, that is the model for the underlying deterministic process $f$.
Interestingly, as shown by \cite{Beek2021}, the predictive variance can also
be decomposed into interpolation error and correction due to noise: $$\vn
(\x) = r(\x) + \pvar \left( c(\x, \x) - \vecc(\x)^\top \CN^{-1} \vecc
(\x) + \sum \limits_{i = 1}^n \vecc(\x)^\top \zeta_i \vecc
(\x)^\top \right),$$ with details on the computation of $\zeta_i$ in the same
article.
In terms of theoretical analysis for the prediction error and convergence,
some existing results with constant noise variance can be found, e.g.,
in \cite{Lederer2019,Garnett2022}. 

So far this assumes that the noise function $r$ is known, which is seldom the
case in practice. One exception is with Monte Carlo error, see e.g., \cite
{Picheny2013a}. Another appears in the context of hyperparameter tuning, where
the accuracy of the output is controlled by an additional parameter. This may
be the size of the data set on which a given machine learning model is
trained, see e.g., \cite{klein2015towards}.

\subsection{Inference for heteroscedastic GP modeling}
\label{ssec:hetGP}

With the variance process $r(\x)$ seldom known, a variety of options have been
proposed to learn it from the data. Increasing slightly in complexity, a
parametric model may be used, such as $r(\x) = \exp(h(\x))$ with $h(\x)$ from a
parametric family, for example, simple polynomials, see e.g., \cite
{Boukouvalas2014} or exponential family \cite{Le2005}. 

For more flexibility, many works entertain a second GP to model the log
variance, ensuring positive variance at the same time, starting from \cite
{Goldberg1998}. The main caveat is that these log variances are unobserved,
thus being latent variables that must be estimated as well. \citet
{Goldberg1998} rely on a fully Bayesian approach and Markov chain Monte Carlo
(MCMC) methods to do so. Later works took approximate inference methods,
expectation-maximization (EM) \cite
{Kersting2007,boukouvalas2009learning}, variational approximation \cite
{Titsias2011} or maximum likelihood \cite{Binois2018}. Even if the learning
is simplified by the underlying smoothing effect of using a GP on the log
noise, the number of latent variables grows with the number of observations,
unless sparse GPs are used, see e.g., 
\cite{snelson2006variable}. They rely on inducing points to summarize the
 entire distribution. 

A perhaps simpler idea is to learn $r(\x)$ from empirical estimates of the
variance, which can only be computed if replication is present in the data.
Denote $y_i^{(j)}, 1 \leq j \leq a_i$,  $a_i$ repeated observations at a
given $\x_i$, with mean $\bar{y_i}$, then $\hat{\sigma}^2_i = \frac{1}
{a_i - 1} \sum \limits_{j = 1}^{a_i} \left(y_i^{(j)} - \bar{y_i}\right)^2$. The GP on the
log noise variance can then be trained directly on these variance estimates,
with the so-called stochastic kriging model \cite{Ankenman2010}. With
replication, among the total $N$ observations, there are only $n$ unique
designs $\xu_i$, $N = \sum \limits_{i=1}^n a_i$. An additional benefit is
that the predictive equations \eqref{eq:hetGP} are equivalent, see e.g., \cite
{Ankenman2010,Picheny2013a,Binois2018}, to:
\begin{align}
m_n(\x) &= \vecc(\x)^\top (\Cn + \Lan \An^{-1})^{-1} \vecYu,\\ 
s_n^2(\x) &= r(\x) + \pvar \left( c(\x, \x) - \vecc(\x)^\top (\Cn + \Lan \An^{-1})^{-1} \vecc(\x) \right)
\end{align} where $\vecYu = (\bar{y_1}, \dots, \bar{y_n})$, $\Lan =  \sigma^{-2} Diag(r
 (\xu_1), \dots, r(\xu_n)$ and $\An = Diag(a_1, \dots, a_n)$. A direct
 benefit of these expressions is that the computational complexity is much
 reduced, depending on $n$ rather than $N$. In terms of error analysis,
 results are available for stochastic kriging, e.g., by \cite
 {wang2019controlling}.

In between the latent variance and empirical approaches, hybrids have been
proposed to leverage the reduced computational complexity coming with
replicates while allowing fewer (even no) replicates per unique design, see
e.g., \cite{boukouvalas2009learning,Binois2018}. This also raises the question
of how to allocate replicates better than fixing the same amount at every
location, which will be discussed in Section \ref{sec:seqdesign}. 

\subsection{Illustration}

\begin{figure}[htpb]
\includegraphics[width=0.5\textwidth, trim= 0 40 0 15, clip]{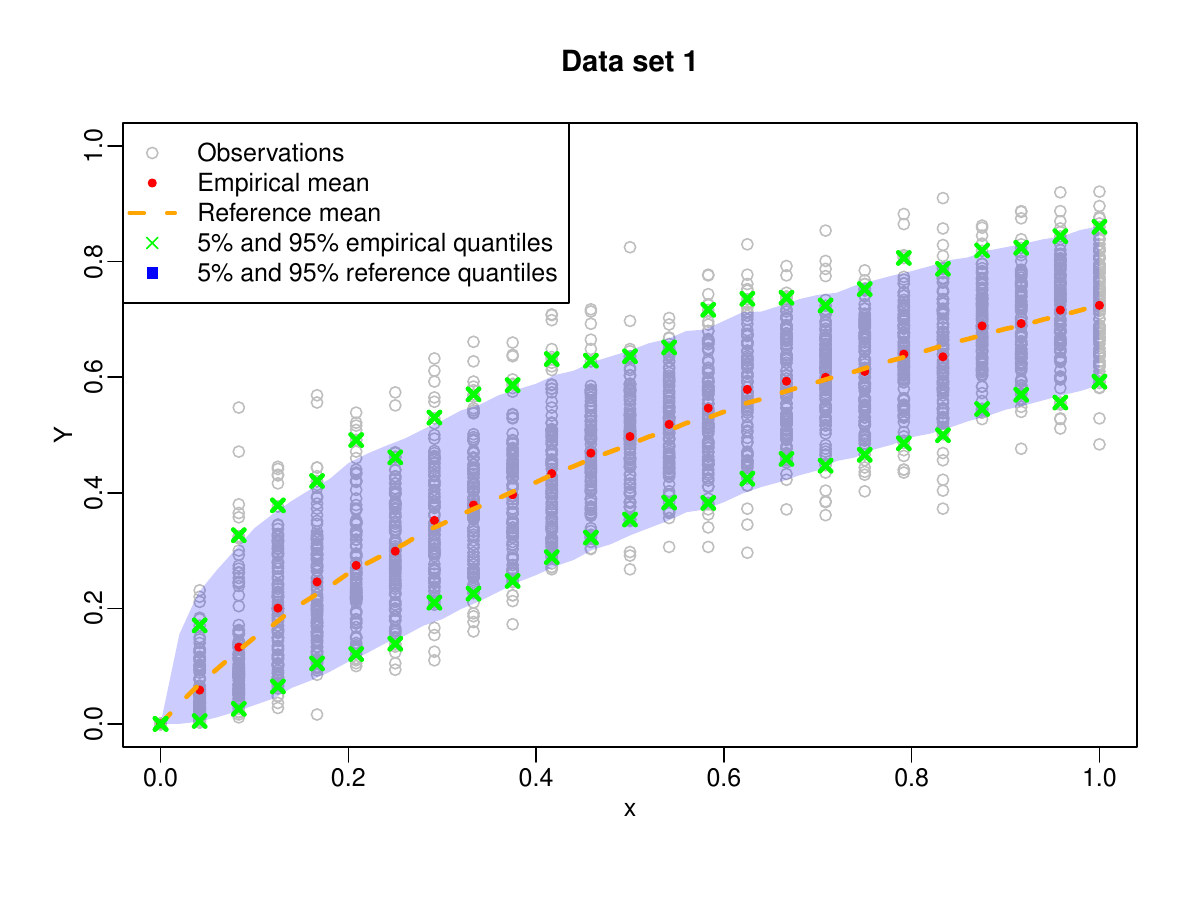}%
\includegraphics[width=0.5\textwidth, trim= 0 40 0 15, clip]{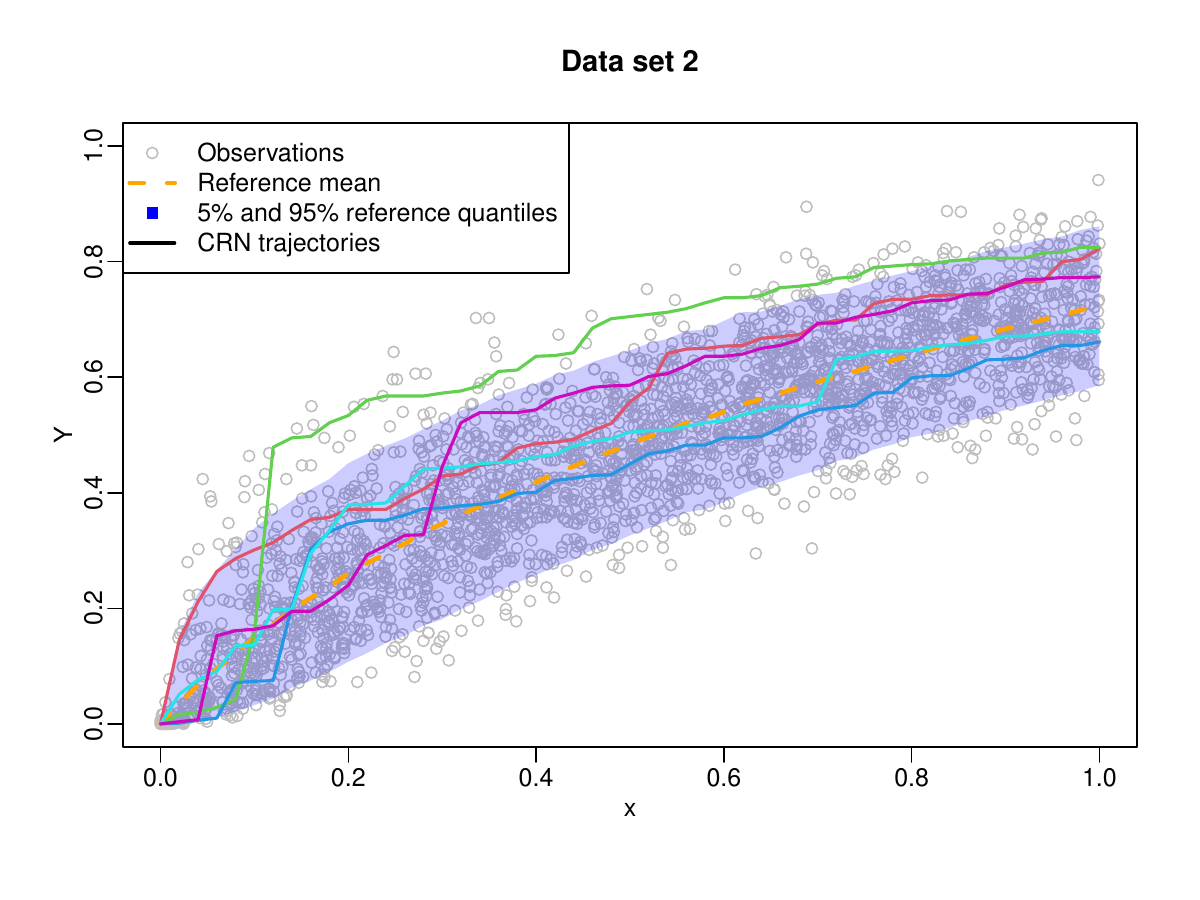}\\
\includegraphics[width=0.5\textwidth, trim= 0 40 0 40, clip]{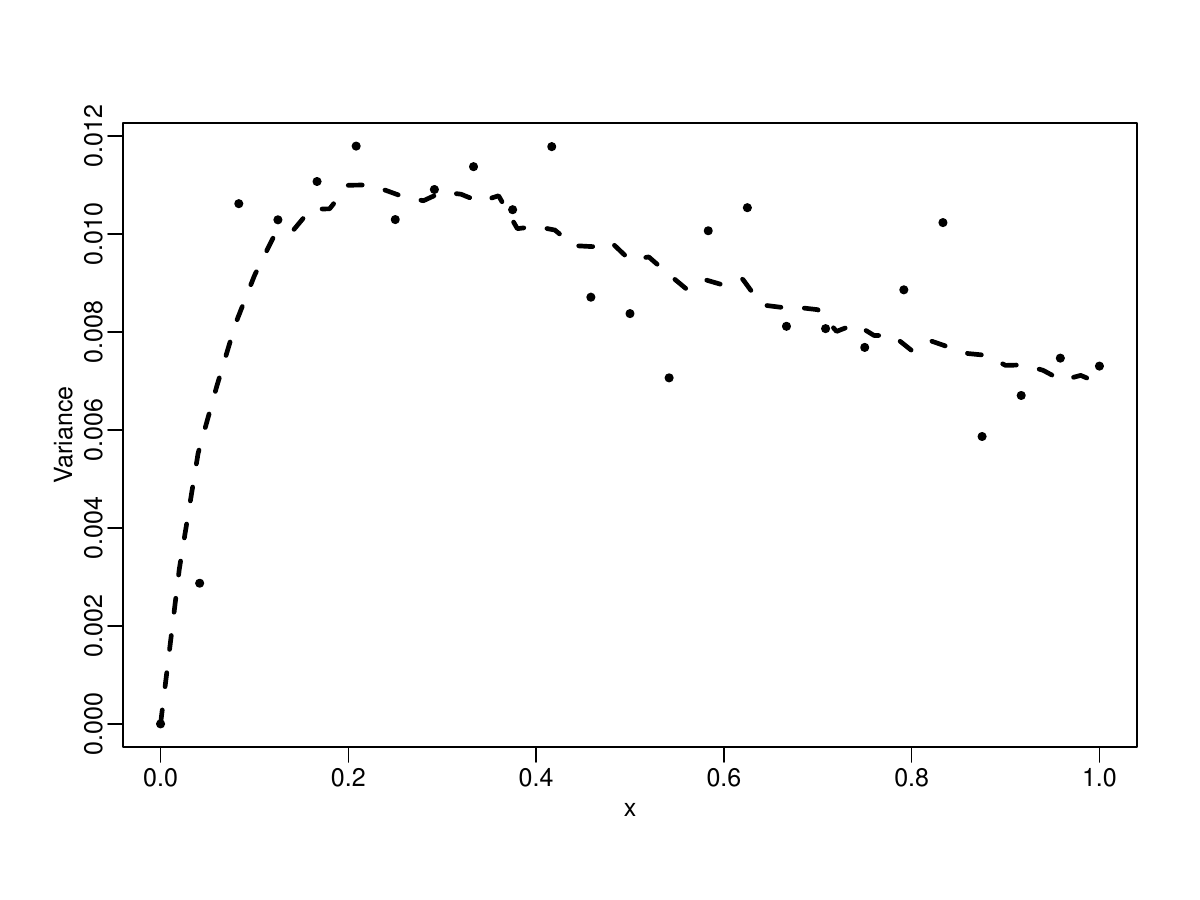}
\includegraphics[width=0.5\textwidth, trim= 0 40 0 40, clip]{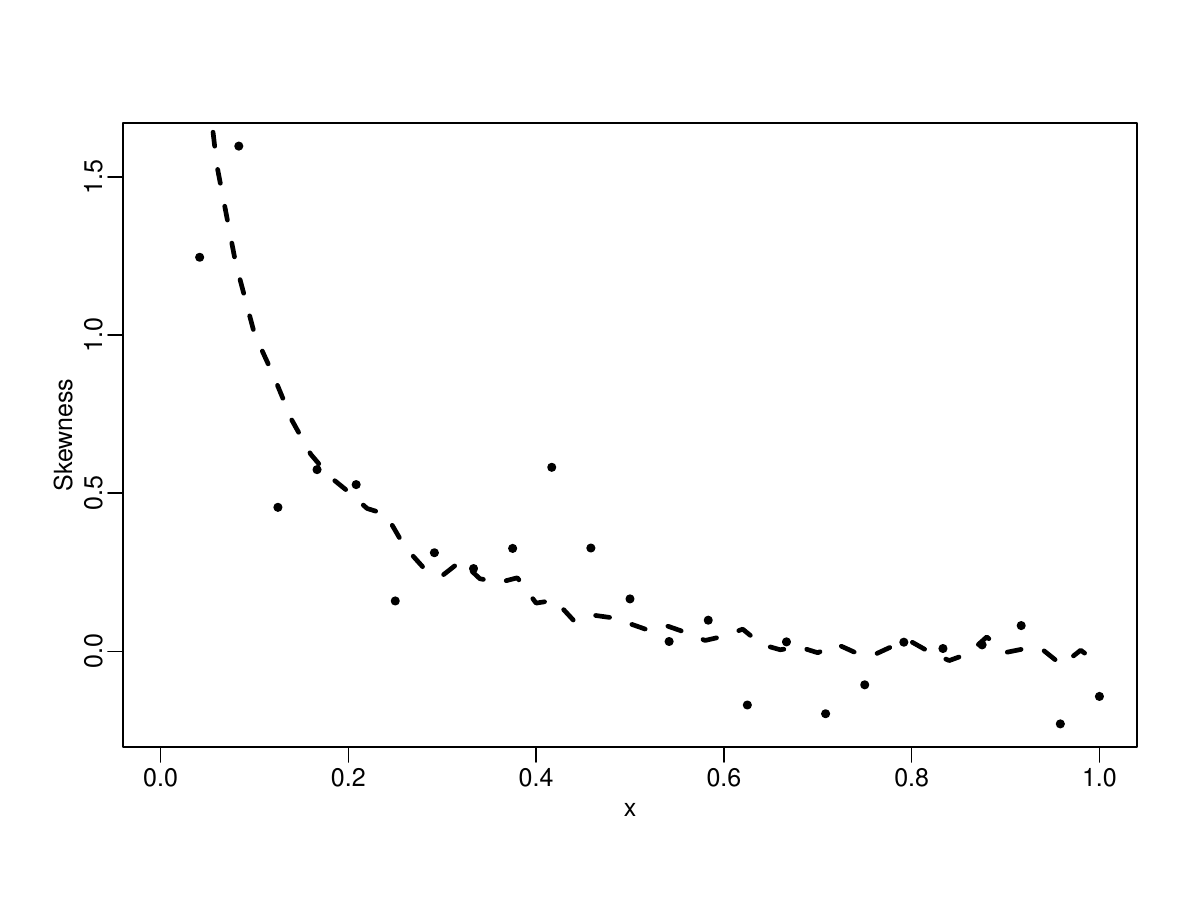}
\caption{SIR simulator data sets. Top left: 2500 observations with
 replications, allowing to compute the estimated mean and quantiles. Top right:
 2500 observations without replications. Five realizations with
 the same seed value across the x-axis are also depicted. Bottom: the empirical
 (points) and reference (dashed lines) variance and skewness for the data set 1.}
\label{fig:testpb}
\end{figure}

We provide an illustrative example of a 1d SIR (Susceptible--Infected--Recovered) problem in Figure \ref
{fig:testpb}. For one data set, we use 25 unique designs each with 100
replicates while there are 2500 unique locations for the second data set. We
use 10000 replicates at 51 unique locations for computing a reference.
Already from this simple data set, it appears that the empirical estimation
is not precise, even using 100 replicates. This even shows for the
reference on the variance and skewness estimation. Without replication,
estimating the variance or higher moments requires learning latent variables.
Trajectories corresponding to the CRN setup, that is fixing the seed values,
are also added for illustration of the difference with the iid setup.

\begin{figure}[htpb]
\includegraphics[width=0.5\textwidth, trim= 0 40 0 15, clip]{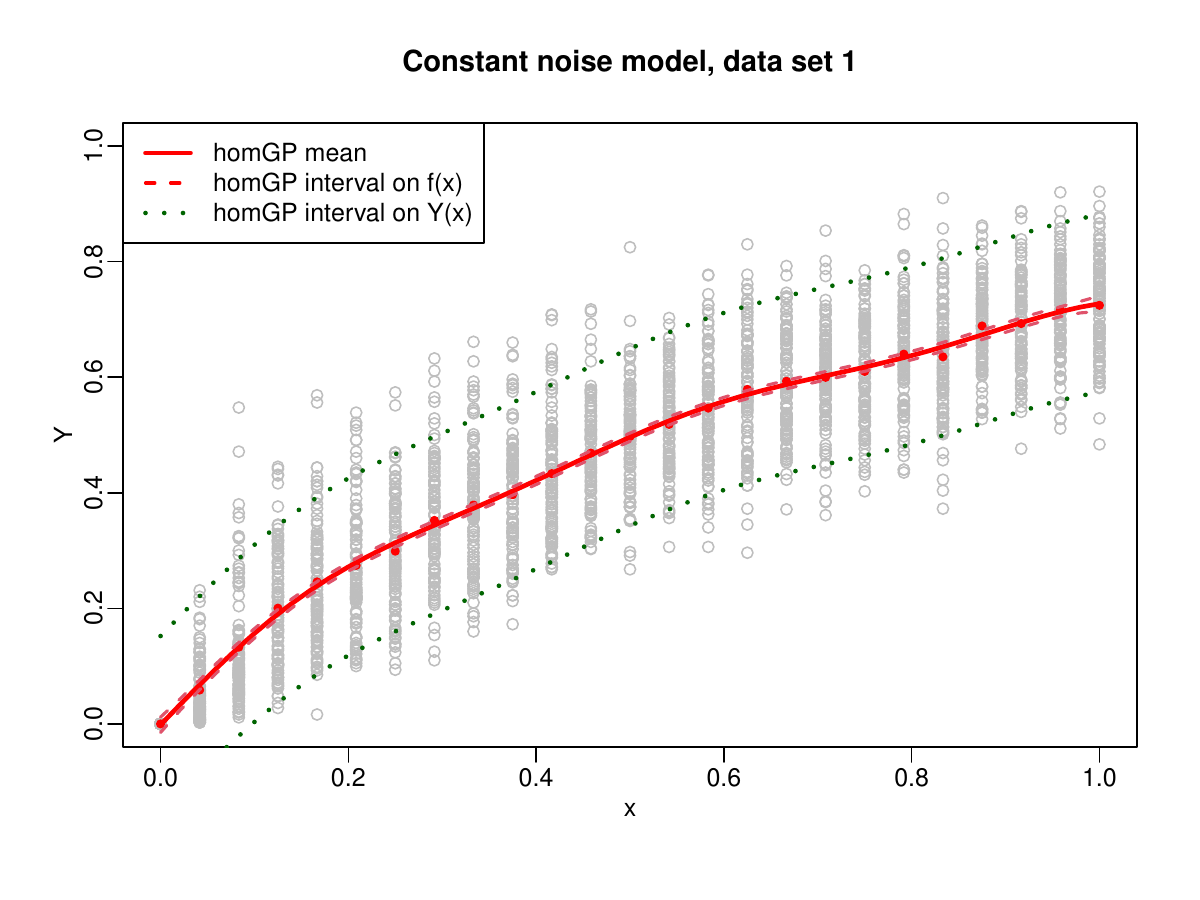}%
\includegraphics[width=0.5\textwidth, trim= 0 40 0 15, clip]{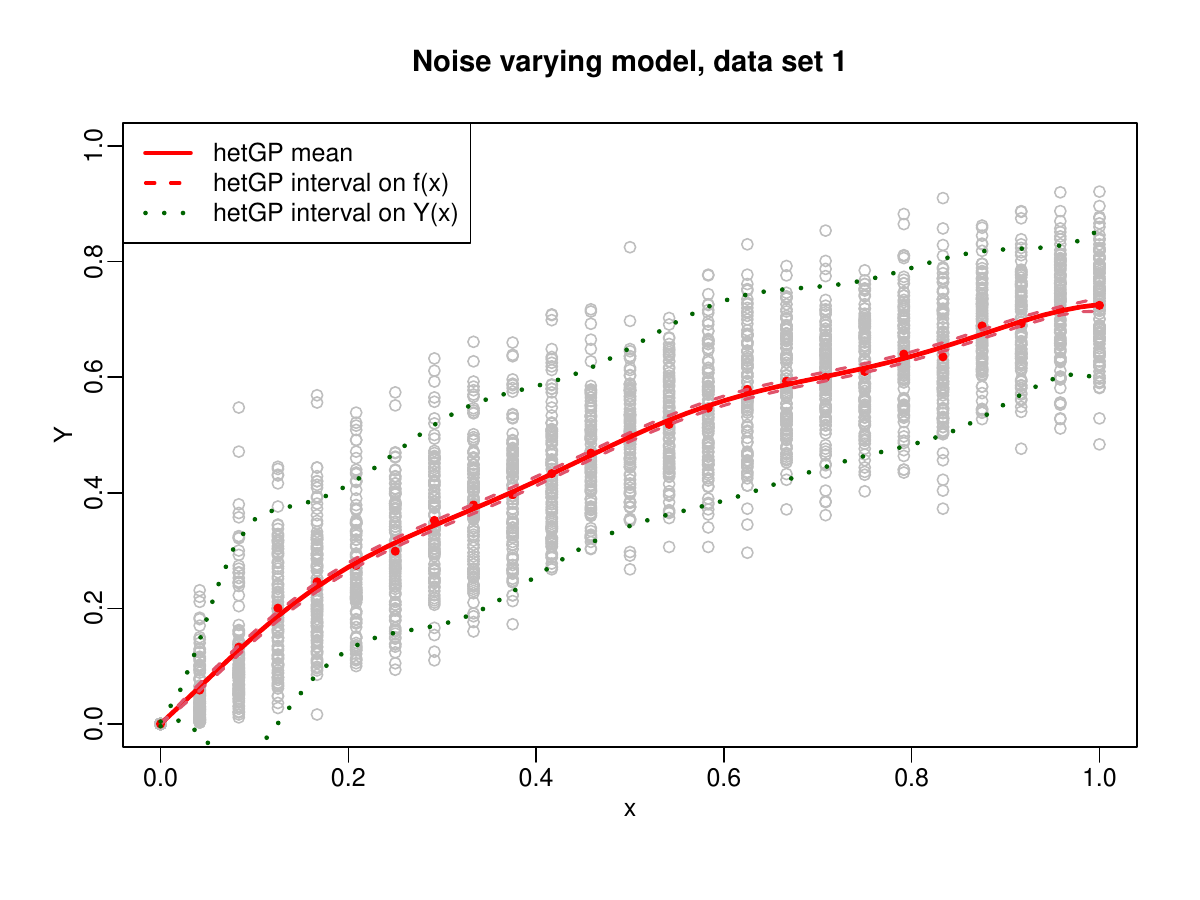}\\
\includegraphics[width=0.5\textwidth, trim= 0 40 0 15, clip]{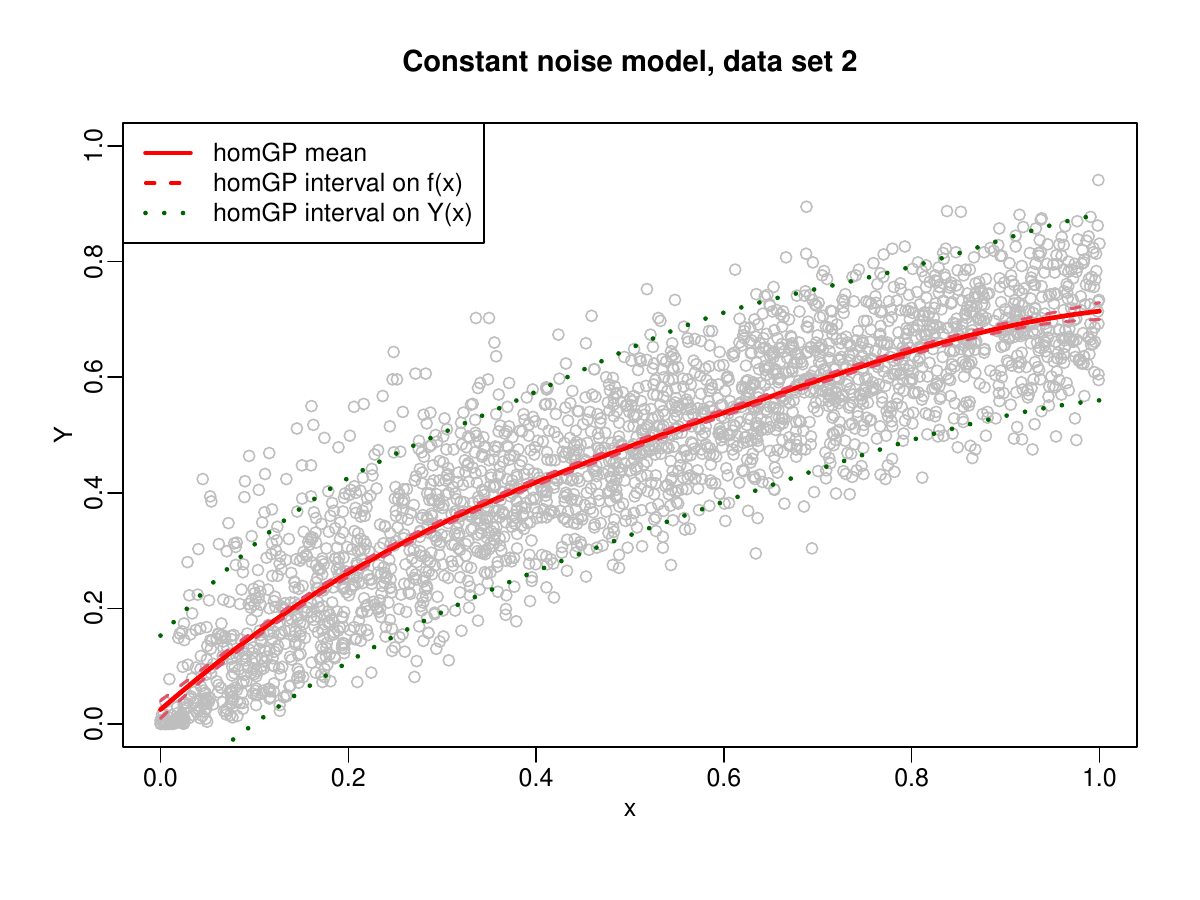}%
\includegraphics[width=0.5\textwidth, trim= 0 40 0 15, clip]{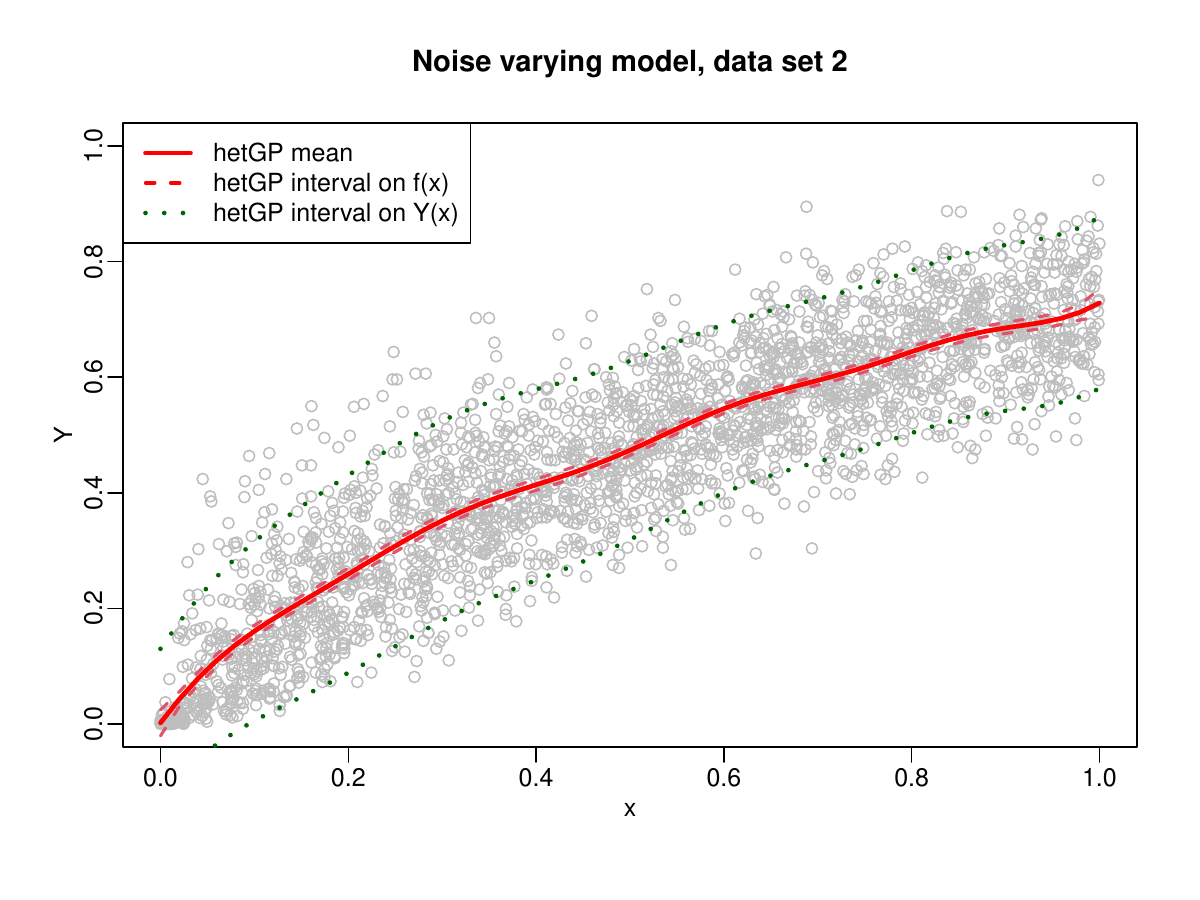}\\
\includegraphics[width=0.5\textwidth, trim= 0 40 0 15, clip]{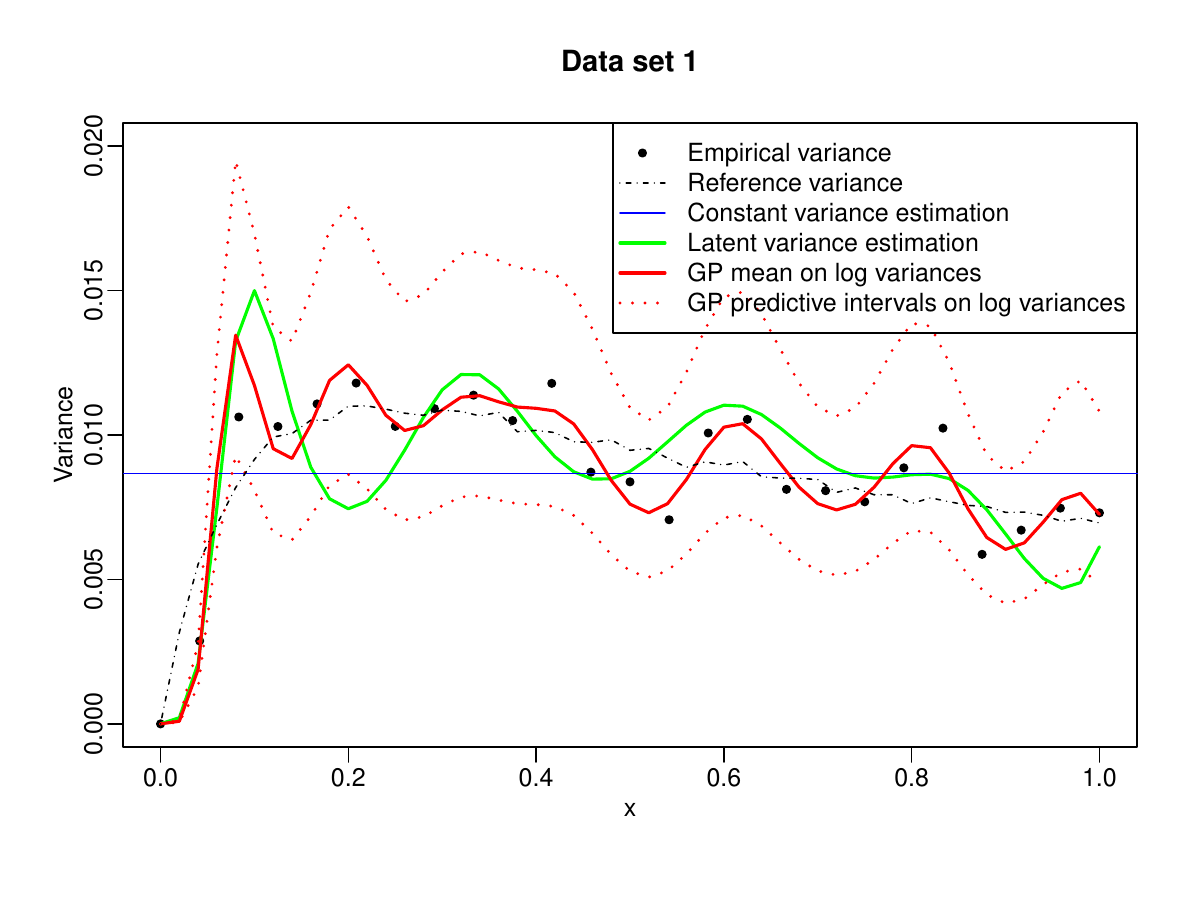}%
\includegraphics[width=0.5\textwidth, trim= 0 40 0 15, clip]{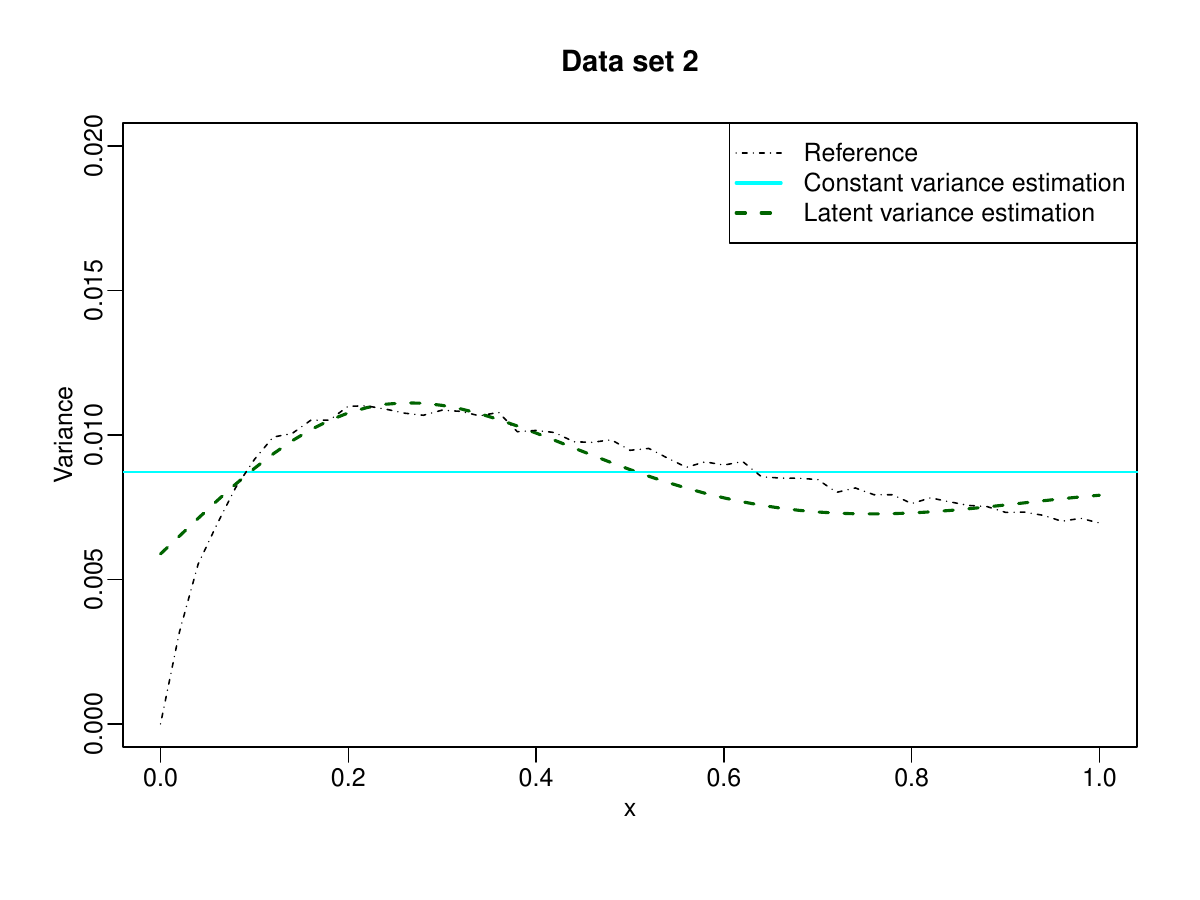}%
\caption{Illustration on a 1d SIR simulator. Homoscedastic versus
 heteroscedastic models predictions on the two data sets. Bottom:
 corresponding variance predictions.}
\label{fig:hgps}
\end{figure}

Then we fit the first GP models on these two data sets, using the \texttt
{hetGP} \cite{Binois2021} \texttt{R} \cite{RCore} package, with results in
Figure \ref{fig:hgps}. The predictive means on the 5 trajectories are also
shown. From this relatively simple example, one can observe the ability of a
heteroscedasic model to better represent the black-box at hand, in terms of the
first two moments. As for the 90\% predictive intervals, since the predictions remain
Gaussian, such a model cannot capture higher moments accurately. As will be
illustrated throughout the entire chapter, capturing finer features requires
more data. Even a GP with constant noise is sufficient to estimate the mean
response, then heteroscedastic GP (hetGP) models also captures the noise variance, but even more
complex models are required to go beyond. The difference between the outcomes
for data sets 1 and 2 are slight, but the training times are very different:
0.08 vs 20.2s for the homoscedastic model, 0.093 vs 396s for the
heteroscedastic one. A stochastic kriging model would need fitting two GPs,
and is only possible on data set 1. From the estimation of variances, a
benefit of the stochastic kriging model is the availability of an estimate of the
predictive variance on the log variance. This is also available depending on
the inference procedure with hetGP models. The use of replicates is
beneficial in terms of variance estimation, it better captures the initial
values around $x \approx 0$ while hetGP models are in the right ballpark.

\section{Extensions beyond the Gaussian case}
\label{sec:advnoisy}

For more complex noise structures, more complex models are needed, at the cost
of less analytical tractability. Besides transforming the outputs with a
parameterized monotonous function as proposed e.g., by \cite
{de1997bayesian,Snelson2004}, we detail some of the variations. A discussion
from the RKHS regularization viewpoint is proposed e.g., by \cite
{Aravkin2015} to improve robustness to outliers. 

\subsection{Some options}

\subsubsection{Student-t variations and non-Gaussian likelihoods}

A first modulation of the Gaussian noise assumption is with Student-t noise,
with larger tails. Using Student-t processes, \cite
{shah2014student,wang2017extended}, analytical expressions may be obtained
but can revert back to GP ones when the number of observations grows,
limiting the interest to smaller data sets.

A more popular assumption is to use a Gaussian process coupled with Student-t
noise, see e.g., \cite{vanhatalo2009gaussian,hartmann2019laplace,Lyu2021}. Nevertheless, the
likelihood is no-longer Gaussian and the posterior distribution cannot be
obtained in closed form. It can still be approximated, using for instance the
Laplace approximation of the posterior, as detailed e.g., by \cite
{Rasmussen2006}. This later approximation is also used with Poisson
distribution for count data, e.g., by \cite{flaxman2015fast}. Another typical
use if for classification with Gaussian processes, with the Bernoulli
distribution.

A wider range of likelihood functions is possible when using MCMC for
approximating the posterior function, made computationally feasible by
relying on variational inducing points, as proposed, e.g., by \cite
{flaxman2015fast}. Coupling variational inducing points and variational
inference for non-Gaussian likelihoods is the approach proposed e.g.,
by \cite{hensman2015mcmc,hensman2015scalable} and available for instance in
the \texttt{GPflow} \cite{Matthews2017}, \texttt{GPyTorch} \cite
{gardner2018gpytorch} or \texttt{GPJax} \cite{Pinder2022} \texttt
{Python} packages. The point of the variational framework is to approximate the intractable posterior distribution by a simpler parameterized one. These
parameters are jointly optimized with the other model parameters by minimizing the Kullback-Leibler divergence between approximate and true posterior via the so-called evidence lower bound (ELBO). Combined with the latent variables idea from Section \ref
{sec:noisygps}, one can model non-Gaussian heteroscedastic models, as
proposed by \cite{saul2016chained} with Student-t or even beta distributed
noises. In terms of flexibility, the generalized lambda distribution, with
its four parameters, can also model a variety of non-Gaussian distributions
and be cast within the same variational framework. This distribution has also
been used with polynomial chaos expansions latent functions by \cite
{zhu2021emulation}. They offer a peek at the underlying phenomenon with the
ability to learn the latent variables.

It remains that these recent generalization capabilities are enabled by
several layers of approximations. Taking a step back, if the goal is really
to learn the tails of unknown noise distributions, then perhaps looking
directly at quantiles (or expectiles) is sufficient, and it as been proposed,
e.g., by \cite{picheny2022bayesian} with two latent GPs.

\subsubsection{Quantile GP modeling}

Instead of assuming i.i.d. additive noise on $\varepsilon_i$, \cite{Plumlee2014} proposed reconstructing replicate variability by estimating quantiles as a function of the inputs. This approach relies on two key principles: (1) with a sufficient number of quantiles, the empirical CDF can be represented as a mixture of point masses, allowing for fast sampling; and (2) quantiles are assumed to be continuous, making them suitable for modeling with a Gaussian Process (GP). The idea is to fit a GP directly to the empirical quantiles, accounting for estimation error, without assuming Gaussian noise. Alternatively, the training data can be modified by estimating multiple quantiles from noisy simulations and augmenting the corresponding quantile levels into the input. Applying a standard GP to this modified data provides similar kriging estimates while avoiding distributional assumptions on the noise. The same formulations in \eqref{eq:hetGP} can be modified as:
\begin{align} \label{eq:quantileGP}
    Y_i(\x_i, \alpha) = f(\x_i, \alpha) + \varepsilon_i, \;\;\; \varepsilon_i \sim N(0, \sigma^2), \;\;1 \le i \le N, 
\end{align}
where, $0<\alpha<1$ is the quantile level. For instance, \cite{Picheny2013a} also considered the GP predictive quantiles, which are available in closed form. If the level of the quantile is not decided
beforehand, then it can be an additional parameter to the approach, as entertained by \cite{fadikar2018calibrating}. Taking a loss function tailored for quantiles and latent input noise
variables, \cite{Quadrianto2009} propose a quantile regression framework based on GPs. It is further extended by \cite{picheny2022bayesian} relying on a variational framework. 

\subsubsection{Deep GPs}
\label{sec:deep}

Noisy simulations can alternatively be viewed through the lens of non-stationarity. 
By applying non-linear transformations to a GP, we can create a stochastic process 
with input-dependent noise variance, albeit without explicit control over the noise 
structure. We focus on warped GPs, where input and/or output are transformed to 
model complex data. Warping functions, such as the monotonic neural-net sum of $\tanh$ 
functions in \cite{Snelson2004}, add flexibility while maintaining the full
probabilistic framework for inference. This model can be further improved by placing
a non-parametric GP prior on the warping function itself, as described in 
\citep{lazaro2012bayesian}, where the inference is performed using variational methods.
Deep GPs are generalizations of the warped GPs. A deep GP is composed of multiple hidden layers of latent variables, where each node serves as the output for the layer above and the input for the layer below. The observed outputs are positioned at the bottom layer, while Gaussian processes control the mappings between layers. A single-layer deep GP (also known as GP-Latent Variable Model or GP-LVM) was first proposed by \cite{lawrence2003gaussian} and a hierarchical version of GP-LVM was introduced
in \cite{lawrence2007gplvm}. \cite{damianou2013deep} later developed the deep GP framework 
and the variational framework to fully infer the model parameters for any number of layers. 
Currently, a number of \texttt{Python} and \texttt{R} packages exist for fitting a deep GP 
model, like for instance \texttt{GPyTorch} \cite{gardner2018gpytorch}, \texttt{GPflux} \cite{dutordoir2021gpflux} or \texttt{deepgp} \cite{Booth2024}.

\subsection{Illustration}

\begin{figure}[htb]
\includegraphics[width=0.5\textwidth, trim= 0 30 0 15, clip]{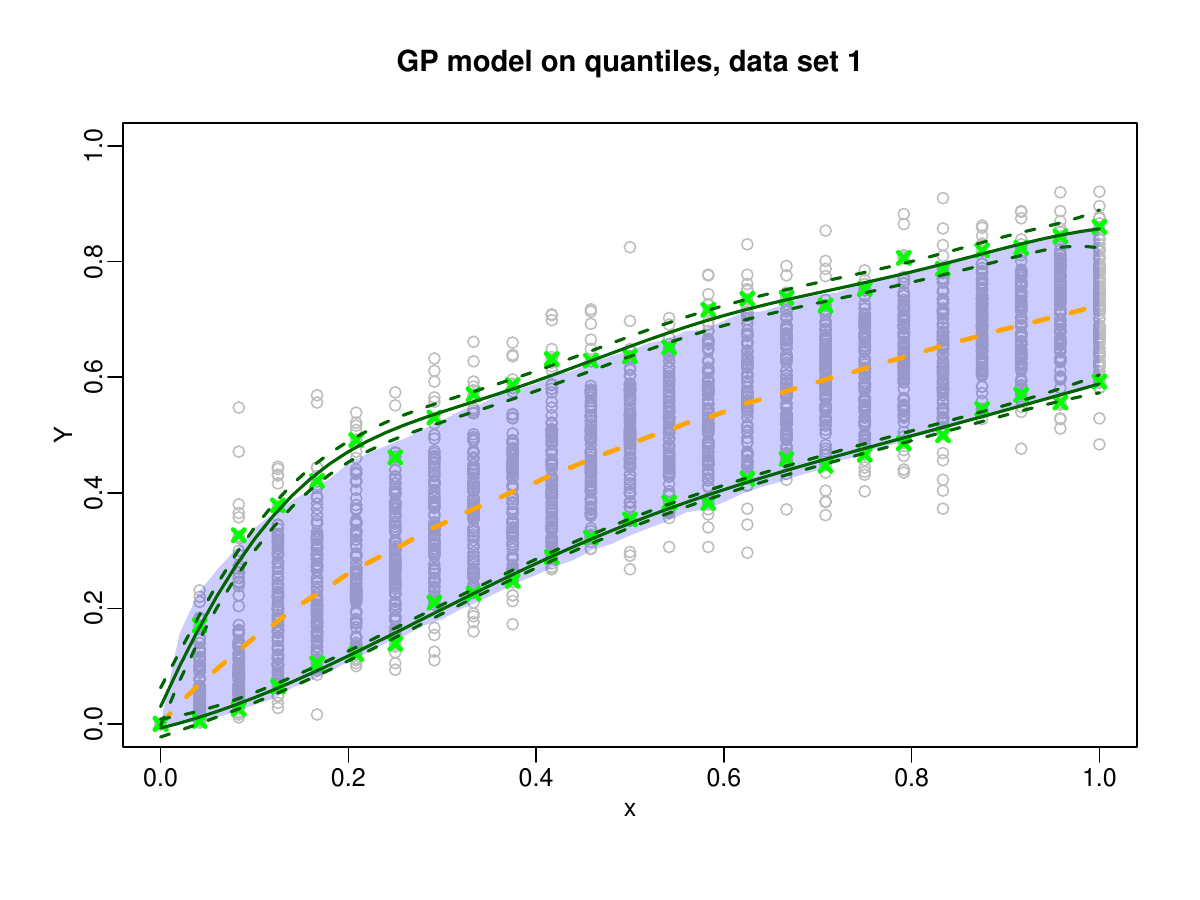}%
\includegraphics[width=0.5\textwidth, trim= 0 30 0 15, clip]{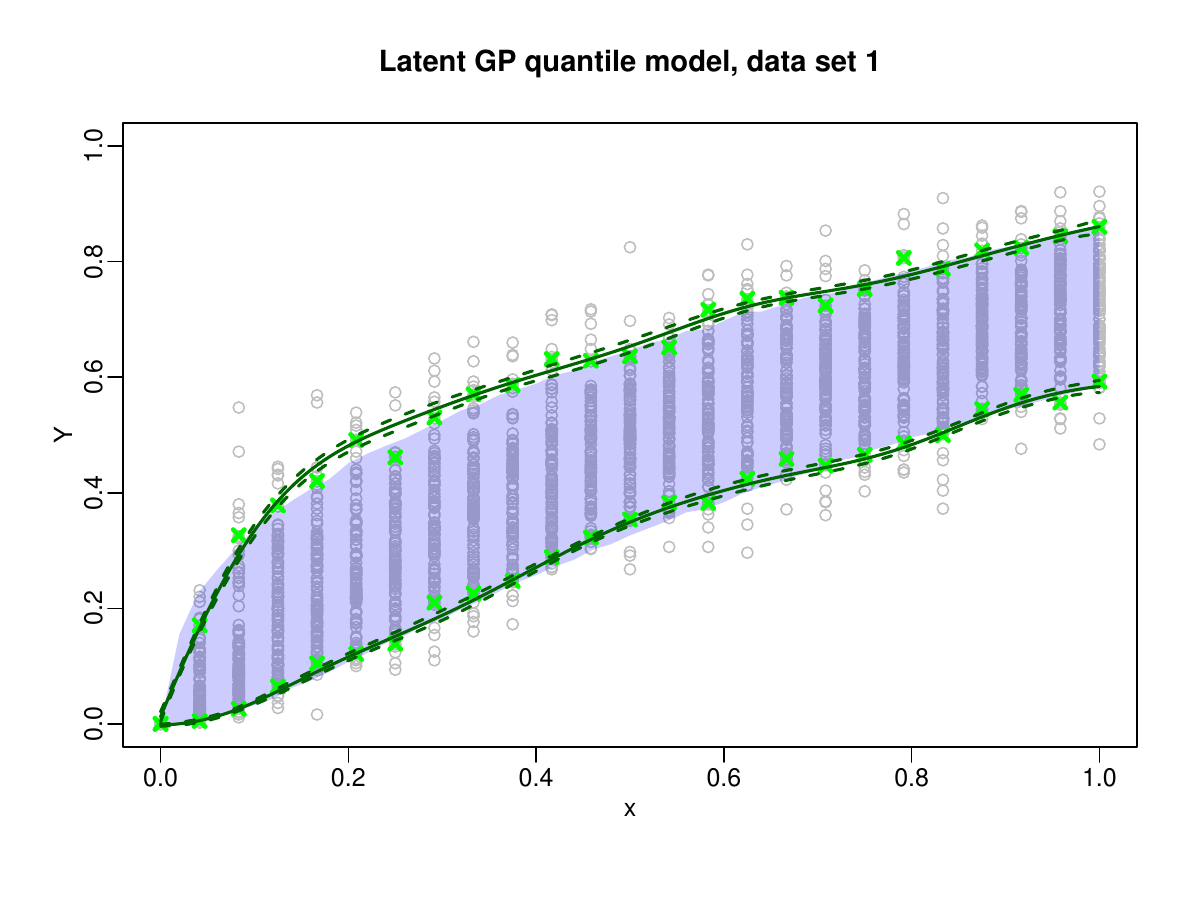}\\
\includegraphics[width=0.5\textwidth, trim=0 30 0 15, clip]{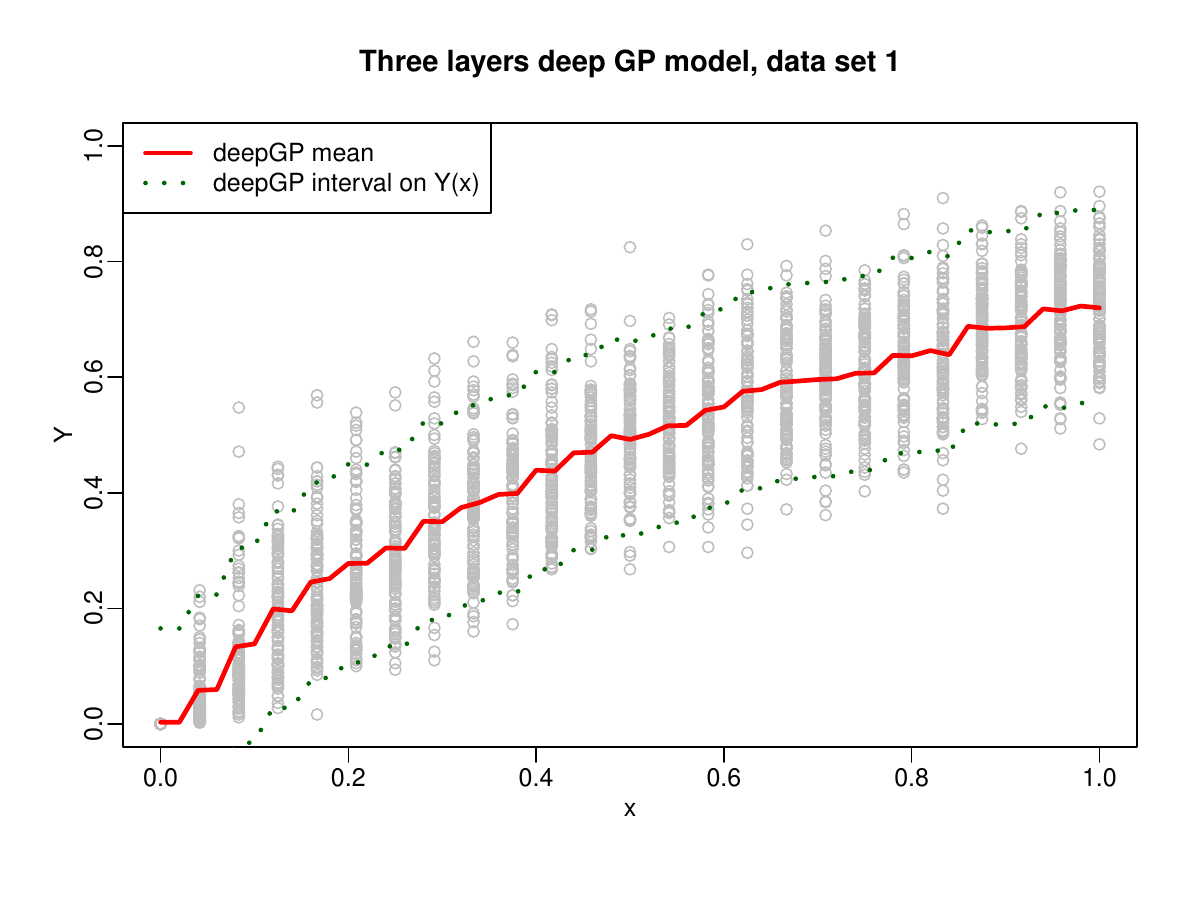}%
\includegraphics[width=0.5\textwidth, trim= 0 30 0 15, clip]{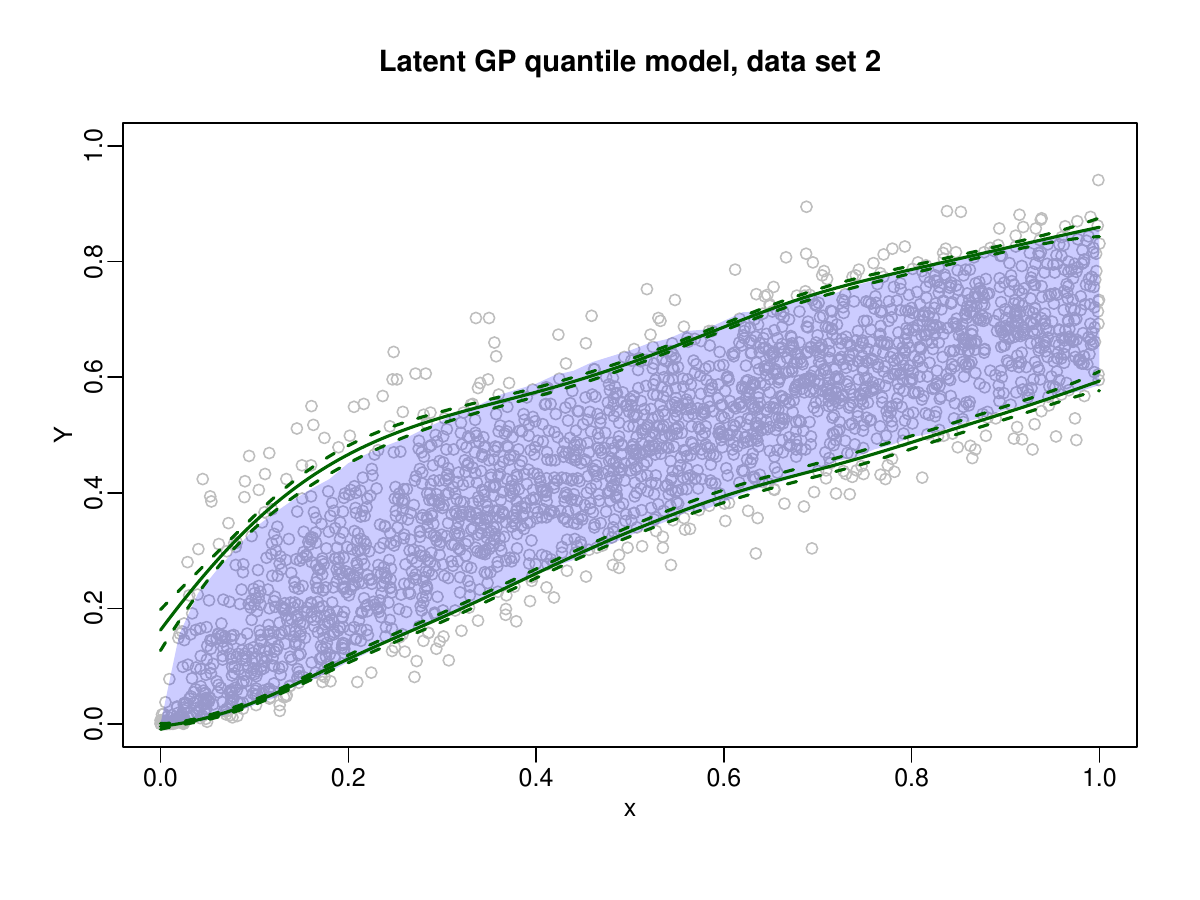}
\caption{Quantile GP models: directly on empirical quantiles (top-left) and using
 latent variables (right). Reference quantiles are in shaded blue. Deep GP model on data set 1 (bottom-left).}
\label{fig:qgps}
\end{figure}

Here we focus on directly modeling the quantiles of the data set, relying either on replicates or latent variables. The outcomes are represented in
Figure \ref{fig:qgps}. For the latent quantile GP, we rely on \cite
{picheny2022bayesian}. On data set 1, the results of both approaches are similar, with broader prediction intervals for the GPs on empirical quantiles. Without replication in data set 2, only the latent quantile model can be applied, which manages to fit the actual quantiles well, except at the origin. Replicates with no or small variance can force the model to take it into account. We fit three layers deep GP models to both data sets using \texttt{deepgp} package in \texttt{R}. However, the model predictions are not significantly different between the two data sets. Figure~\ref{fig:qgps} shows the model fit on the data set with replication.

After presenting broad options of modeling capabilities, the specific choice
for a given application would take into account various aspects. First, the
amount of available data will drive the choice, the most complex models
requiring more data. Using an appropriate noise model is also important to
accurately represent the data, especially if this noise structure is known. Finally, the structure of the design of
experiments may also restrict options: with replication, estimating moments
is facilitated for modeling them directly. Otherwise, inducing points and variational inference of latent quantities may be necessary to keep a reasonable computational cost.

Next we will review the options to improve the model for specific targets.

\section{Adapted sequential design strategies}
\label{sec:seqdesign}

While noise modeling brings more complexity, it also introduces a
regularization component which may actually alleviate numerical issues
(covariance matrix conditioning) and sequential design procedures
(the acquisition function value is smoother as it is not zero at evaluated
designs).

Before delving into specific methods for various targets, some general
techniques are possible, see for instance \cite{Jalali2016, rojas2020survey}. A first
naive one is simply to ignore the noise and do exactly as in the noiseless
case. This is probably fine as long as the signal to noise ratio is high.
Increasing this latter is achievable by fixing a number of replicates for each
evaluation. A more careful approach is to filter out the noise, by using the
prediction by the noisy GP model as inputs for a deterministic GP employed in
the sequential approach. This is called reinterpolation, see, e.g., \cite
{Forrester2008}. 

We continue by detailing several typical applications: global accuracy and
optimization. Estimating a probability of failure, see e.g., \cite
{hao2021novel}, or calibration \cite{baker2022analyzing}, or a level set, see
e.g., \cite{Lyu2021}, are other related tasks, and so is sensitivity
analysis \cite{mazo2021trade}. We refer to e.g., \cite
{kleijnen2015design,Santner2018,Ginsbourger2018} and references therein for
more general aspects. 

\subsection{Global accuracy}

Broadly speaking, getting a globally accurate surrogate model amounts to have
good space-filling properties, see e.g., \cite{Gramacy2020}. There a
reasonable criterion to adopt is the integrated mean square prediction error
(IMSPE): $\alpha_{IMSPE}(\x) = \int_\Xset \vn(\x) d\x$. Optimizing it
sequentially attempts to reduce the predictive variance everywhere, without
putting too much effort on the frontier of the domain. Alternatively,
targeting model estimation is studied for instance by \citet
{Boukouvalas2014} with the Fisher information.

The heteroscedastic noise context brings a twist, to spread the effort between
areas where the noise variance is more or less large. When the input set
$\Xset$ is discrete, the optimal allocation of replicates over designs is
given by \cite{Ankenman2010}. For continuous inputs, \citet{Binois2018a} show
that replication may be preferred to adding new designs depending on $r
(\x)$ and propose further strategies to promote replication. A key idea is to
look ahead, meaning taking into account the effect of future decisions. That
is, reducing the number of alternatives to keep tractability, to assess
whether or not evaluating a new design now and then replicating is better
than waiting and replicate a few times before adding a new design. Since IMSPE only depends on the predictive
variance, which does not depend on the observations, this is a favorable
context for these less myopic options. Batch extensions are detailed e.g., in \cite
{zhang2022batch}.

\subsection{Optimization}

Global accuracy criteria are usually related to the predictive variance or
covariance matrix only. This is less simple for other goal oriented
approaches such as those geared towards optimization -- over the mean in this
section -- that require adaptations. Notable exceptions where the noisy setup is transparent include the probability of improvement, knowledge
gradient or entropy based criteria, albeit these two latter families require
approximations to be computed and an estimate of $r(\x)$. 

The difficulty for
the popular expected improvement (EI) criterion is to define the current best value $T$:  $\alpha_{EI}(\x) = \Esp{\max(0, T -Y(\x)|y_1, \dots,y_N)}$. Indeed, noise realizations may trick the surrogate model, leading to an overly optimistic estimation of this target. It will then require more samples for correction. This effect may be more pronounced for multi-objective optimization, where a set of Pareto optimal solution are needed.
Properly identifying the current best estimate from the surrogate model is called the identification step by \cite{Jalali2016}, which is
required at least at the final iteration to return a solution to the user. 

For EI, a simple
adaptation is to plug-in the predictive mean instead of function values, at
the evaluated designs \cite{vazquez2008global} or over $\Xset$ \cite
{Gramacy2011}. An offset may be added for more robust estimation of this
plugin, combined to a penalty for the noise variance at $\x_{n+1}$ by \cite
{huang2006global}. With additional noisy constraints, a Bayesian version of
EI is proposed by \cite{letham2019constrained}. Thorough discussions available in \cite{frazier2018bayesian,Gramacy2020,Garnett2022}
details relationships between EI, KG or entropy
 infill criteria, and extensions to related topics: non myopic, batch or
 multi-fidelity criteria. Last but not least, \cite{astudillo2019bayesian} insists on
 modeling the actual output directly, rather than modeling a sum, an integral
 or even a composition. 

Criteria dedicated to the noisy context include the expected quantile
improvement \cite{picheny2010noisy}, looking at the improvement over a
quantile. Notice that the upper confidence bound (UCB) criterion can be seen
as directly optimizing a GP quantile for a given level. Benchmark of the
above criteria are presented for instance by \cite
{picheny2013benchmark,Jalali2016}. \citet{Jalali2016} insist on the
importance of a replication strategy, interleaving an additional stage to the
usual procedure and on an appropriate identification strategy to actually return a good solution at the end. 

So far only standard GPs, including the heteroscedastic version have been
considered. In the case of Student-t processes, dedicated criteria have also
been obtained in closed form, for multi- \cite{Herten2016} or
single-objective \cite{xie2017heteroscedastic} optimization. In the context
of CRNs, we refer to \cite{Pearce2022} for exploiting the relationship
between seeds to obtain the best mean results and to \cite
{fadikar2023trajectory} for finding good individual realizations
(irrespective of the mean behavior).

Only a few criteria exist to directly choose the number of replicates when
evaluating, see e.g., \cite{Jalali2016,Beek2021,binois2021portfolio}. But two stage
approaches have been proposed in the operation research literature, first new
points are evaluated, then the remaining budget is used to replicate at the
most promising locations. Some rules have been derived for selecting and
ranking solutions, sometimes under the term optimal computing budget
allocation(OCBA) \cite{chen2000simulation,quan2013simulation}, or for
multi-objective optimization counterparts, e.g., \cite
{Hunter2017,rojas2020survey,gonzalez2020multiobjective}.

Another way to see this is to consider the cost associated with
a noise variance reduction, hence cost-aware strategies can be employed, see
e.g., \cite{Garnett2022}. We mention another simple strategy when the signal
to noise ratio is low: enforcing a minimal variance reduction for a new
evaluation. This helps ensuring that there will be sufficient information to
actually modify the prediction significantly, hence reducing the time for
BO.

\subsection{Robust optimization}

A related topic when dealing with uncertainty is to edge against it by the
formulation of the problem: rather than considering simply the expected
behavior, a threshold on the variance may be used (or estimating a
mean-variance Pareto front), a worst case approach may be adopted, optimizing a
given quantile or taking a chance constrained framework. Nevertheless these
are generally tackled considering environmental variables (necessary, e.g., for
worst-case frameworks). 

Using noisy surrogates, optimization on quantiles is tackled e.g., by \cite
{browne2016stochastic,picheny2022bayesian}, with a model of the quantiles and
an EI-like infill criterion, or by \cite{makarova2021risk} for a weighted sum
of the mean plus variance. Additional criteria looking at the
(conditional) value at risk are studied for instance by \cite
{chen2016efficient} within a stochastic kriging framework, when \cite
{picheny2022bayesian} uses latent GPs for directly modeling the quantiles
with entropy and Thompson sampling acquisition functions. Finally, a
multi-objective version is described, e.g., by \cite{rivier2022surrogate}.

\subsection{Illustration}

We show two strategies: one selecting individual evaluation at every
iteration, but targeting replicates, one choosing simultaneously the next
design with its replication budget, corresponding to a variance reduction of at least 10\%. Which one is preferable depends on the level of
concurrency for evaluations, the signal to noise ratio or the time for single
evaluations. 

\begin{figure}[htb]
\includegraphics[width=0.5\textwidth, trim= 0 30 0 15, clip]{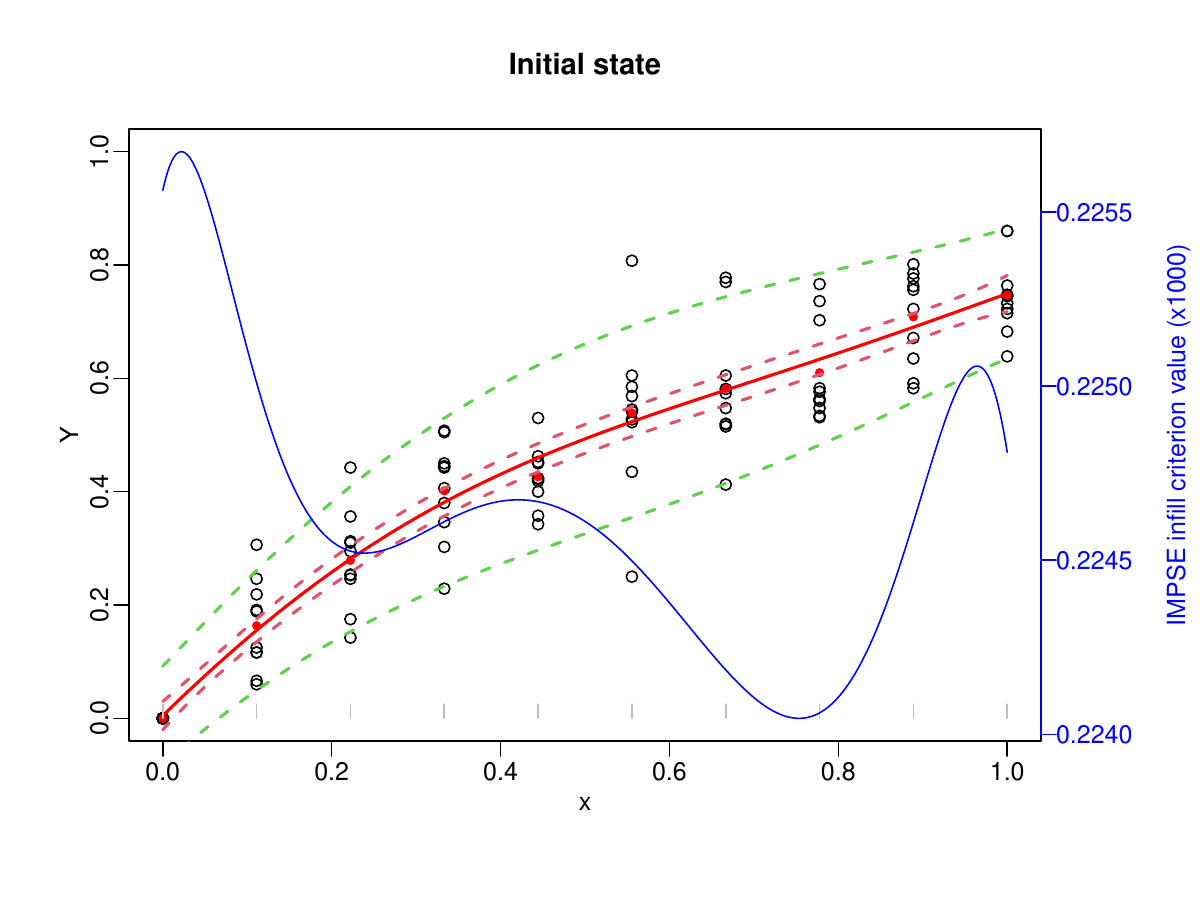}%
\includegraphics[width=0.5\textwidth, trim= 0 30 0 15, clip]{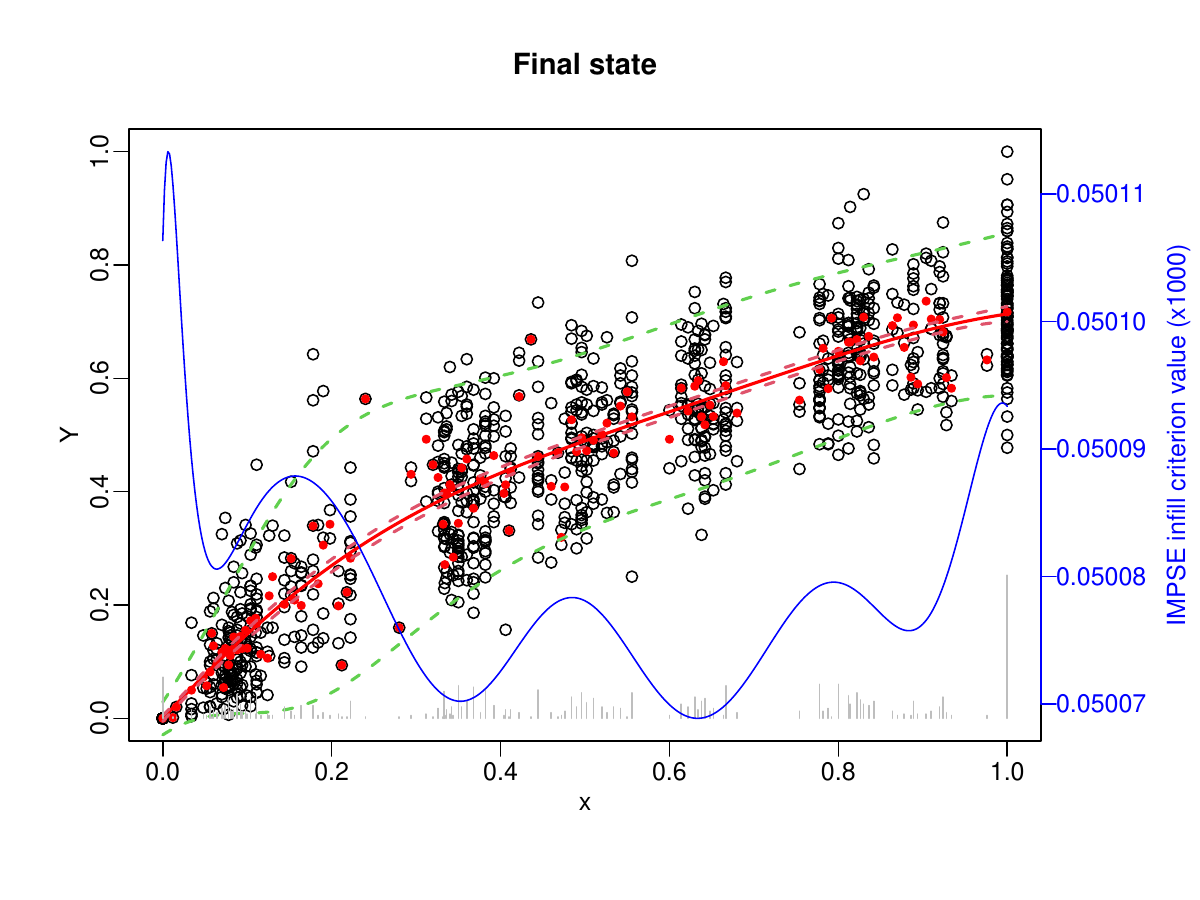}\\
\includegraphics[width=0.5\textwidth, trim= 0 30 0 15, clip]{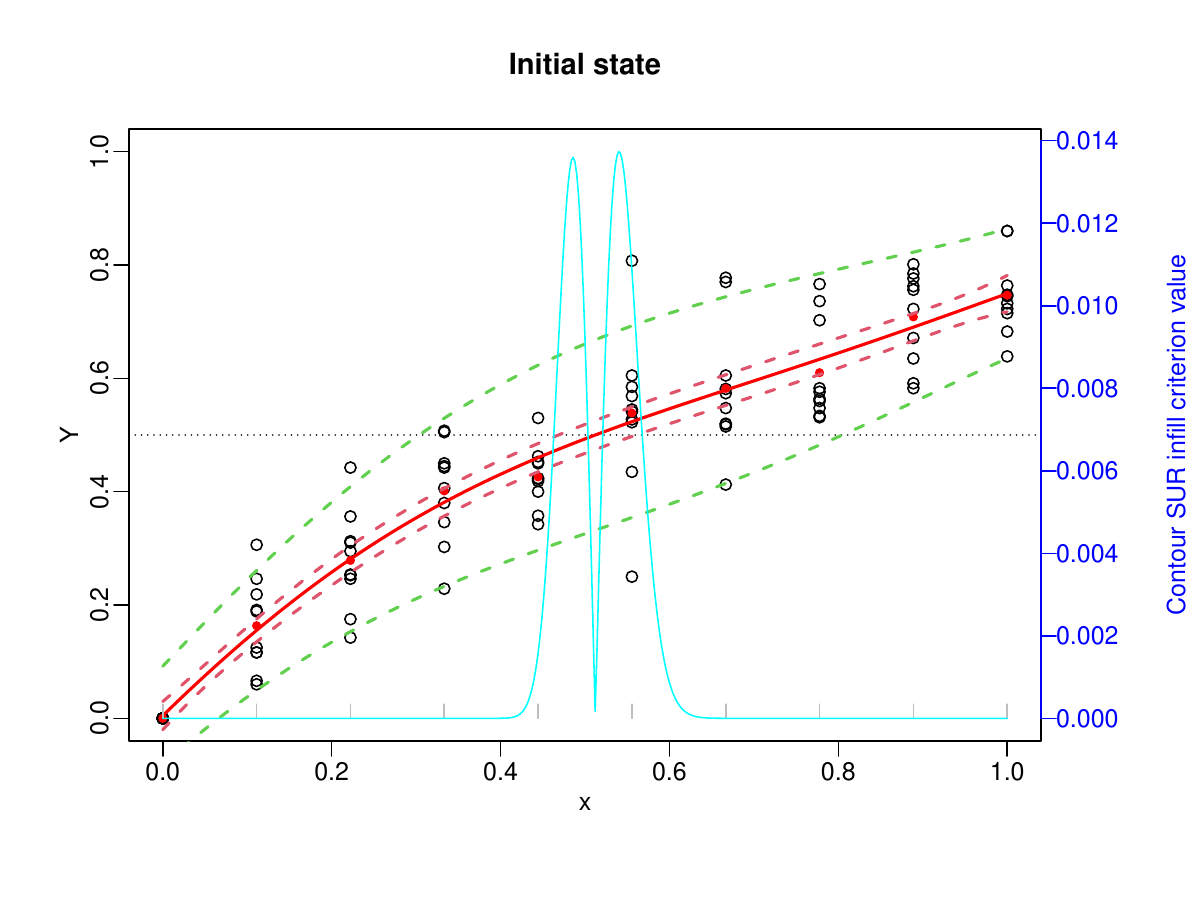}%
\includegraphics[width=0.5\textwidth, trim= 0 30 0 15, clip]{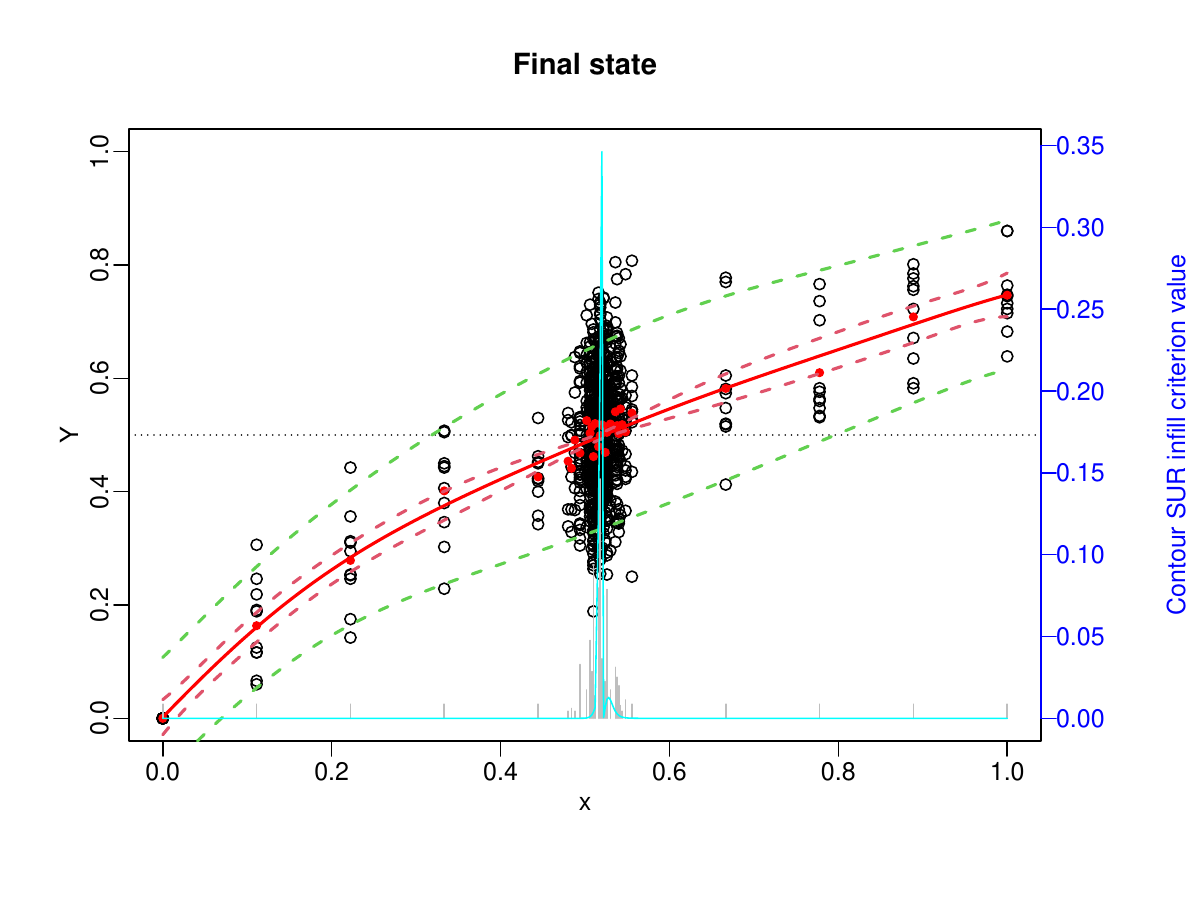}%
\caption{Sequential learning illustration. Left: initial state, right: state
 after using IMSPE (global accuracy, top) with empirical means in red when there are replicates. Bottom: state for 0.5 level set
 estimation. Gray vertical bars represent the amount of replication per design.}
\label{fig:levels}
\end{figure}

A first task is to reach global accuracy, hence we follow the approach relying
on looking ahead with IMSPE to lightly enforce replication.  For the second task, we consider a level set exercise to identify where $f=0.5$. The contour stepwise uncertainty reduction (SUR)  criterion proposed by \cite{Lyu2021} is used. Similar to the knowledge gradient idea, it focuses on the effect of a new observation on the estimation of the contour.
Results for these two setups are presented in
Figure \ref{fig:levels}.
For IMSPE, a total of 900
points are selected sequentially, one per iteration, resulting in $\approx
120$ unique designs. With the contour SUR, only $\approx 30$ iterations are necessary to reach
a total of 1000 evaluations. Both succeed in improving the initial estimate
towards the goal, either placing replicates in a more or less space-filling
way for IMSPE, except on the left where the noise variance is smaller while
for the level set the replicates are on both sides of the crossing.

\section{Conclusion and perspectives}
\label{sec:concl}

In this chapter, we reviewed a variety of methods for modeling noisy
simulators. Some benefit from or require replication, while others make no such
assumption. The former may be simpler to interpret and re-implement, while the
latter are available with advanced implementations and inference frameworks.
Both frameworks have been adapted for sequential designs.

We focused mostly on GP models, but other random fields model may be relevant as
well, see e.g., \cite{kozubowski2013multivariate}. Similar in structure to
their neural networks counterparts, deep GPs use Gaussian processes for their
layers, see e.g.,  \cite
{damianou2013deep,dunlop2018deep,sauer2023vecchia}, and have been adapted to
sequential learning as well \cite
{hebbal2021bayesian,yazdi2022fast,sauer2023active} but, as discussed in Section \ref{sec:deep}, without dedicated
treatment for the noise. Further studies of their performance in the noisy setup would be needed, perhaps with additional focus on the inference procedure. 

In terms of remaining road blocks, most are shared with regular GPs, such as
coping with non-stationarity or high-dimensional inputs. More specific ones
in the noisy context include model selection depending on the amount of
available data or more fine-tuned handling of replicates.

\end{document}